\documentclass[12pt]{article}
\usepackage{amsmath, amsthm}
\usepackage{graphicx, color}

\newbox\struttbox\setbox\struttbox=\vbox to 0.75\baselineskip{}
\def\strutt{\copy\struttbox}
\def\pagetable#1\endpagetable{\begin{table}[p]
    \vbox to 1\textheight{\centering#1}\end{table}\clearpage}

\begin{document}
\title
{Hypohamiltonian planar cubic graphs with girth five}
\author{
Brendan D. McKay\\
{\small Research School of Computer Science}\\[-0.4ex]
{\small Australian National University}\\[-0.4ex]
{\small ACT 2601, Australia}\\
{\small \tt E-mail: bdm@cs.anu.edu.au}
}
\date{}

\maketitle

\begin{abstract}
A graph is called hypohamiltonian if it is not hamiltonian but becomes
hamiltonian if any vertex is removed.
Many hypohamiltonian planar cubic graphs have been
found, starting with constructions of Thomassen in 1981.
However, all the examples found until now had 4-cycles.
In this note we present the first examples of hypohamiltonian planar cubic
graphs with cyclic connectivity five, and thus girth five.
We show by computer search that the smallest members of this class
are three graphs with 76 vertices.
\end{abstract}

{\small {\bf Keywords:} graph generation; hypohamiltonian graph; planar graph;
  cubic graph}

\smallskip

{\small \textbf{MSC 2010.} 05C10, 05C30, 05C38, 05C45, 05C85.}

\section{Introduction}

A graph is called \textit{hypohamiltonian} if it is not hamiltonian but becomes
hamiltonian if any vertex is removed.
Such graphs exist on orders 10 (the Petersen graph), 13, 15, 16, and all 
orders from 18 onwards~\cite{AMW}.
It is elementary to show that hypohamiltonian graphs must be 3-connected
and have cyclic connectivity at least~4.

A substantial body of literature is devoted to finding hypohamiltonian
graphs with various properties.  The first planar hypohamiltonian graph,
with 105 vertices, was found by Thomassen~\cite{Thom76}; the present
smallest known order is 40 vertices but smaller orders have not been
ruled out~\cite{Joo14}.  

Planar hypohamiltonian graphs can be cubic, as first shown by Thomassen
with an example on 94 vertices~\cite{Thom81}.  The smallest examples
found so far have 70 vertices; the first example by Araya and
Wiener~\cite{Araya11} and six more by Jooyandeh
and McKay~\cite{Joo13}.

Planar hypohamiltonian graphs can also have girth 5, and in this case
the smallest order is proven to be 45~\cite{Joo14}.  However, it was
not known until now whether planar hypohamiltonian graphs can be
both cubic and of girth 5.  Our purpose is to present the first examples.

\section{The results}

Using the program~\texttt{plantri}~\cite{Brink05,Brink07} we generated all
planar cubic graphs of girth 5 and cyclic connectivity at least 4 on up
to 76 vertices.  No theory is known which could efficiently restrict the 
search to a small subclass sure to contain the hypohamiltonian graphs, so
we just tested all of them.

The computational task was daunting as more than
$10^{13}$ graphs are involved and 76 vertices is large enough that
a primitive backtrack search for hamiltonian cycles takes up to 
several hours per graph.  Fortunately we have a program (unpublished)
specifically designed for sub-cubic graphs that takes only 10
microseconds on average to find a hamiltonian cycle if there is one,
and 15 milliseconds on average to rigorously prove that there is none.
Even then the task took about 8 years of cpu time.

\begin{figure}
\centering
\includegraphics[scale=0.4]{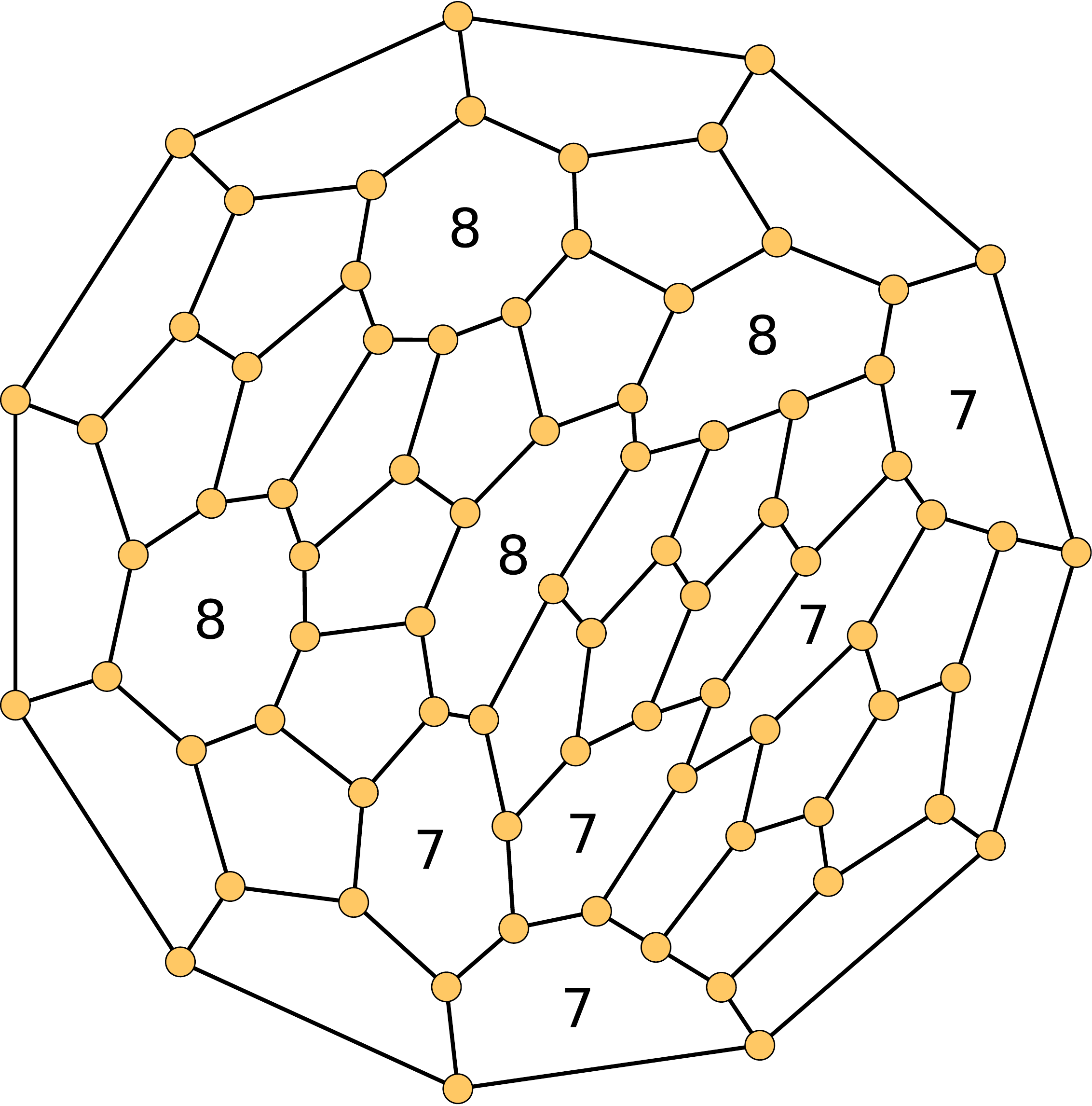}
\caption{Hypohamiltonian graph with face count $5^{30}\,7^5\,8^4\,11$}
\label{fig1}
\end{figure}

\begin{figure}
\centering
\includegraphics[scale=0.4]{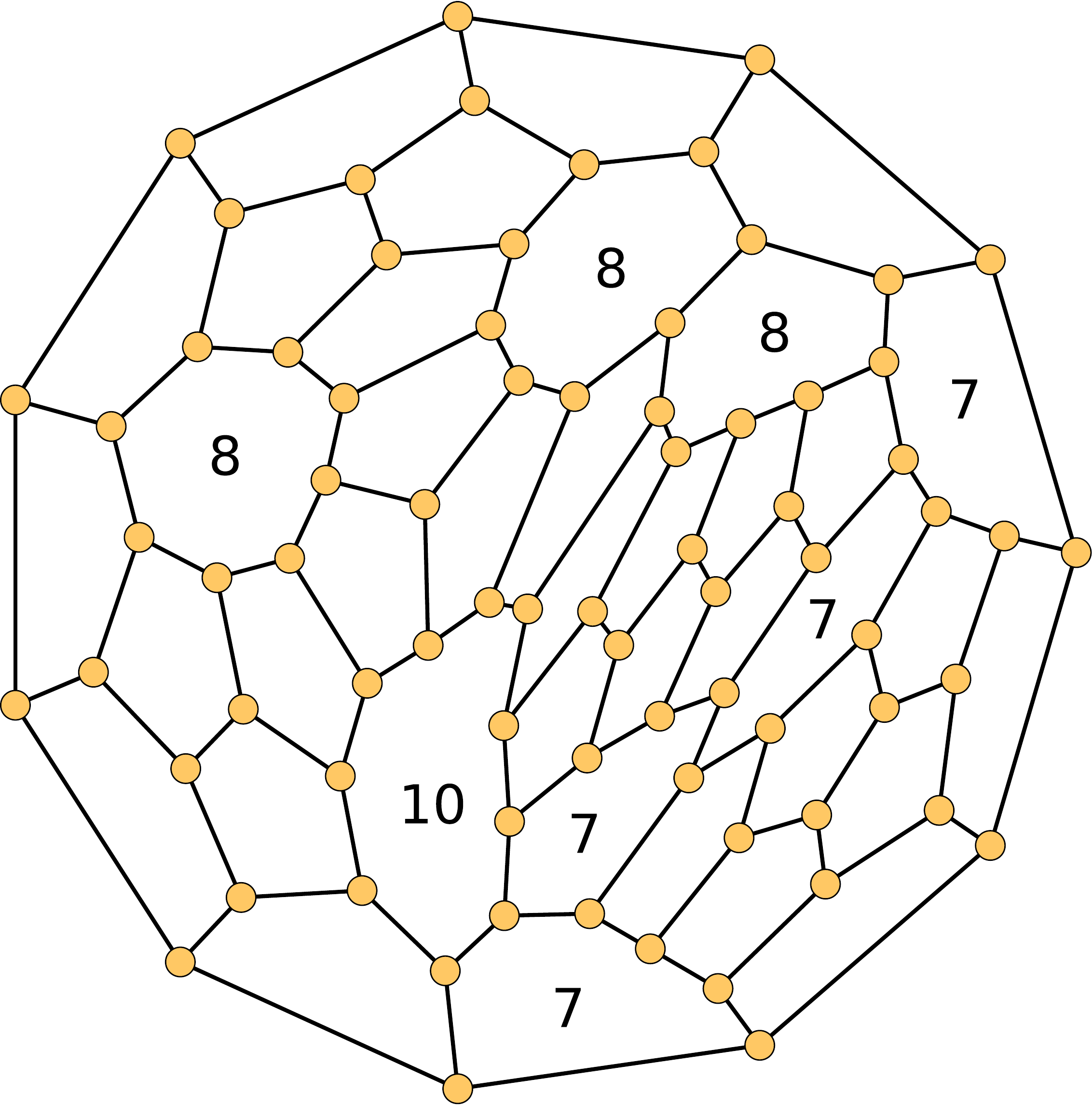}
\caption{Hypohamiltonian graph with face count $5^{31}\,7^4\,8^3\,10\,11$}
\label{fig2}
\end{figure}

\begin{figure}
\centering
\includegraphics[scale=0.4]{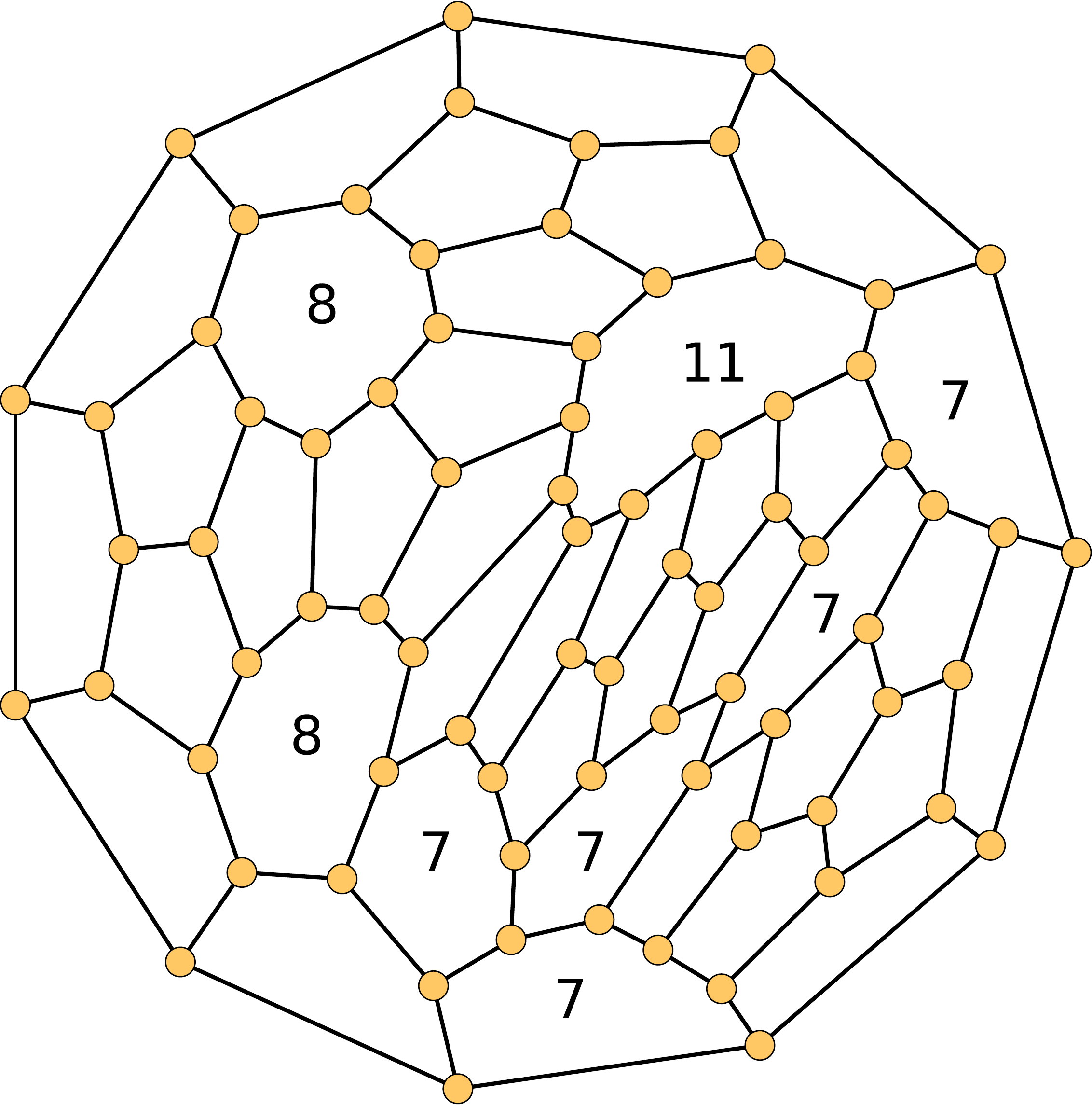}
\caption{Hypohamiltonian graph with face count $5^{31}\,7^5\,8^2\,11^2$}
\label{fig3}
\end{figure}

The results are shown in Table~\ref{tab}.  For at most 74 vertices
there are no hypohamiltonian graphs, but for 76 vertices there
are three.  All of them have cyclic connectivity 5.  The
graphs themselves are shown in Figures~\ref{fig1}--\ref{fig3}, 
which were drawn by CaGe~\cite{CaGe}.

Note that the graphs have different face counts despite having some
similarities. 
The non-hamiltonicity of these graphs does not follow immediately
from Grinberg's condition, so we tested them with several independent
programs to remove the possibility of error.
Contrary to our usual
experience of extremal graphs, none of them have any non-trivial
automorphisms. 
We did not find a way to generalize these examples
to larger sizes despite trying several approaches, so the problem 
of finding an infinite family remains open.

\pagetable\vss
\setlength{\tabcolsep}{4mm}
\begin{tabular}{c|ccccc}
 $n$ & $C_4(n)$ & $N_4(n)$ & $C_5(n)$ & $N_5(n)$ & $H(n)$ \\[0.2ex]
\hline\strutt
 20 & 1 & & 1 & & \\
 22 & 0 & & 0 & & \\
 24 & 1 & & 1 & & \\
 26 & 1 & & 1 & & \\
 28 & 3 & & 3 & & \\
 30 & 4 & & 4 & & \\
 32 & 12 & & 12 & & \\
 34 & 23 & & 23 & & \\
 36 & 2 & & 71 & & \\
 38 & 4 & & 187 & & \\
 40 & 22 & & 627 & & \\
 42 & 84 & & 1970 & & \\
 44 & 376 & 1 & 6833 & 1 & \\
 46 & 1579 & 3 & 23384 & 1 & \\
 48 & 6751 & 1 & 82625 & 0 & \\
 50 & 27969 & 3 & 292164 & 3 & \\
 52 & 115423 & 6 & 1045329 & 6 & \\
 54 & 467948 & 12 & 3750277 & 2 & \\
 56 & 1882184 & 49 & 13532724 & 22 & \\
 58 & 7496828 & 126 & 48977625 & 37 & \\
 60 & 29667311 & 214 & 177919099 & 31 & \\
 62 & 116710547 & 659 & 648145255 & 194 & \\
 64 & 457122502 & 1467 & 2368046117 & 298 & \\
 66 & 1783850057 & 3247 & 8674199554 & 306 & \\
 68 & 6941579864 & 9187& 31854078139 & 1538 & \\
 70 & 26950926431 & 22069& 117252592450 & 2566 & \\
 72 & 104455609591 & 50514 & 432576302286 & 3091 & \\
 74 & 404298188921 & 137787 & 1599320144703 & 13487 & \\
 76 & 1563255455769 & 339804 & 5925181102878 & 22274 & 3
\end{tabular}
\begin{align*}
 n &= \,\text{the number of vertices} \\
 C_4(n) &= \,\text{the number with cyclic connectivity exactly 4}\\
 N_4(n) &= \,\text{the number of those which are not hamiltonian}\\
  C_5(n) &= \,\text{the number with cyclic connectivity exactly 5}\\
 N_5(n) &= \,\text{the number of those which are not hamiltonian}\\
 H(n) &= \,\text{the number which are hypohamiltonian}
\end{align*}
 The numbers $C_4(n)$ and $C_5(n)$ first appeared in~\cite{Brink05}.
\vss
\caption{\label{tab} Counts of planar cyclically 4-connected cubic graphs of girth 5}
\endpagetable


\begin{thebibliography}{99}


\bibitem{AMW}
R.\,E.\,L. Aldred, B.\,D. McKay, and N.\,C. Wormald, Small hypohamiltonian graphs,
\emph{J. Combin. Math. Combin. Comput.} \textbf{23} (1997) 143--152.

\bibitem{Araya11}
M. Araya and G. Wiener, On cubic planar hypohamiltonian and hypotraceable graphs,
\emph{Electron. J. Combin.} {\bf 18}(1) (2011) \#P85.


\bibitem{CaGe}
G. Brinkmann, O. Delgado Friedrichs, S. Lisken, A. Peeters and N. Van Cleemput,
CaGe---a virtual environment for studying some special classes of plane
graphs---an update,
\emph{MATCH Commun. Math. Comput. Chem.} \textbf{63} (2010) 533--552.

\bibitem{Brink05}
G. Brinkmann and B.\,D. McKay,
Construction of planar triangulations
with minimum degree 5, {\it Discrete Math.} {\bf 301} (2005) 147--163.

\bibitem{Brink07}
G. Brinkmann and B.\,D. McKay,
 Fast generation of planar graphs,
  {\it MATCH Commun. Math. Comput. Chem.} {\bf 58} (2007) 323--357.
  
%










\bibitem{Joo13}
M.~Jooyandeh and B.\,D. McKay,  Hypohamiltonian planar graphs,
Web site at 
\texttt{http://cs.anu.edu.au/$\sim$bdm/data/planegraphs.html}.

\bibitem{Joo14}
M. Jooyandeh, B.\,D. McKay, P. \"Osterg{\aa}rd, V. Pettersson and C.\,T. Zamfirescu,
Planar hypohamiltonian graphs on 40 vertices,
submitted (2013). \texttt{arXiv:1302.2698}




%



\bibitem{Thom76}
C. Thomassen, Planar and infinite hypohamiltonian and hypotraceable graphs,
\emph{Discrete Math.} \textbf{14} (1976) 377--389.


\bibitem{Thom81}
C. Thomassen, Planar cubic hypohamiltonian and hypotraceable graphs,
\emph{J. Combin. Theory, Ser.~B} \textbf{30} (1981) 36--44.


%


%



\end{thebibliography}
\end{document}